\documentclass[reqno,10 pt]{amsart}

\usepackage[all]{xy}

\theoremstyle{plain}
\newtheorem{lemma}{Lemma}[section]
\newtheorem{theor}[lemma]{Theorem}

\theoremstyle{definition}

\theoremstyle{remark}
\newtheorem{rem}[lemma]{Remark}

\newcommand{\Pic}{{\rm Pic\thinspace}}
\newcommand{\bs}{{\rm Bs}}

\newcommand{\vdim}{{\rm vdim\thinspace}}
\newcommand{\edim}{{\rm edim\thinspace}}

\newcommand{\p}{\mathbb{P}}

\newcommand{\co}{{\mathcal O}}
\newcommand{\cl}{{\mathcal L}}

\newcommand{\ci}{{\mathcal I}}

\newcommand{\et}{\tilde{e}}
\newcommand{\hti}{\tilde{h}}

\newcommand{\snvs}{\vspace{-1mm}}

\newcommand{\map}{\rightarrow}

\begin{document}

\title[Very ampleness on blowings-up of $\p^3$]{A note on the very ampleness of complete linear systems on blowings-up of $\p^3$}
\author{Cindy De Volder}
\address{
Department of Pure Mathematics and Computeralgebra,
Galglaan 2, \newline B-9000 Ghent, Belgium}
\email{cdv@cage.ugent.be}
%\thanks{The first author is a Postdoctoral Fellow of
%the Fund for Scientific Research-Flanders (Belgium)
%(F.W.O.-Vlaanderen)}

\author{Antonio Laface}
\address{
Dipartimento di Matematica, Universit\`a degli Studi di Milano,
Via Saldini 50, \newline 20100 Milano, Italy }
\email{antonio.laface@unimi.it}
%\thanks{The second author would like to acknowledge the support of the MIUR of the
%Italian Government in the framework of the National Research
%Project ``Geometry in Algebraic Varieties'' (Cofin 2002)}

%\keywords{Very ampleness, base point freeness, linear systems, fat points, projective space}
%\subjclass{14C20}
\begin{abstract}
In this note we consider the blowing-up $X$ of $\p^3$ along $r$ general points
of the anticanonical divisor of a smooth quadric in $\p^3$.
Given a complete linear system $\cl = |dH - m_1 E_1 - \cdots - m_r E_r|$ on $X$,
with $H$ the pull-back of a plane in $\p^3$ and $E_i$ the exceptional divisor corresponding to $P_i$,
we give necessary and sufficient conditions for the very ampleness (resp. base point freeness and non-speciality) of $\cl$. As a corollary we obtain a sufficient condition for the very ampleness of
such a complete linear system on the blowing-up of $\p^3$ along $r$ general points.
\end{abstract}
\maketitle

\section{Introduction}

In this note we work over an algebraically closed field of characteristic $0$.

Let $P_1, \ldots, P_r$ be general points of the anticanonical divisor of a smooth quadric in $\p^3$
and
choose some integers $m_1 \geq \ldots \geq m_r \geq 0$.
Consider the linear system $\cl'$ of surfaces of degree $d$ in $\p^3$ having
multiplicities at least $m_i$ at $P_i$, for all $i=1,\ldots, r$.
Let $X$ denote the blowing-up of $\p^3$ alongs $P_1, \ldots, P_r$,
and let $\cl$ denote the complete linear system on $X$ corresponding to  $\cl'$.

%The virtual dimension of the system $\cl$ is given by
%\[
%\vdim(\cl):=\binom{d+3}{3}-\sum_{i=1}^r\binom{m_i+2}{3} -1,
%\]
and the system is called {\em special} if $\dim(\cl)>\max\{-1,\vdim(\cl)\}$.

Let $Z$ be a zero-dimensional subscheme of length 2 of $X$,
then $\cl$ separates $Z$ if there exists a divisor $D \in \cl$ such that
$Z \cap D \neq \emptyset$ but $Z \not\subset D$.
The system $\cl$ on $X$ is called very ample if it separates all such $Z$.

The very ampleness of line bundles of blowings-up of varieties has been studied by several
authors, e.g. St\'ephane Chauvin and the first author \cite{CD}, Marc Coppens \cite{MC1,MC2},
Brian Harbourne \cite{BH2}, Mauro C. Beltrametti
and Andrew J. Sommese \cite{BS}.

In theorem~\ref{va} (resp. theorem~\ref{bpf}) we prove that such a system $\cl$ is very ample (resp. base point free) on $X$ if and only if
$m_r > 0$, $d \geq m_1 + m_2+1$ and $4d \geq m_1 + \cdots + m_r +3$ (resp. $d \geq m_1 + m_2$ and $4d \geq m_1 + \cdots + m_r +2$).
A fundamental tool for proving these results is theorem~\ref{ns} which states that a system $\cl$ with
$2d\geq m_1+\cdots+m_4$ is non-special if $d\geq m_1+m_2-1$ and $4d\geq m_1 + \cdots + m_r$.

If $r \leq 8$, the points $P_i$ are in general position on $\p^3$ and the
dimension and base locus of $\cl$ on $X$ can be determined using the results from~\cite{DL1,DL2}.

The techniques used in this note are a generalization of the ones in~\cite{DL1,DL2} and make use of the results about complete linear systems on rational surfaces with irreducible anticanonical divisor (see~\cite{BH1,BH2}).

\section{Preliminaries and notation}

Let $\cl_3(d)$ denote the complete linear system of surfaces of degree $d$ in $\p^3$.
Consider a general quadric $\bar{Q} \in \cl_3(2)$ in $\p^3$ and
let $K_{\bar{Q}}$ denote the canonical class on $\bar{Q}$.
Then we know that $-K_{\bar{Q}}$ is just the linear system on $\bar{Q}$ induced by $\cl_3(2)$,
so we can consider $D_{\bar{Q}} \in -K_{\bar{Q}}$ which is smooth and irreducible.

Let $P_1, \ldots, P_r$ be general points of $D_{\bar{Q}}$ and
choose integers $m_1 \geq \cdots \geq m_r \geq 0$.
By $X_r$ we denote the blowing-up of $\p^3$ along the points $P_1 , \ldots, P_r$,
$E_0$ denotes the pullback of a plane in $\p^3$,
by $E_i$ ($i=1,\ldots,r$) we mean the exceptional divisor on $X_r$ corresponding to $P_i$
and $\pi : X_r \map \p^3$ denotes the projection map.

On $\p^3$, we let $\cl_3(d ; m_1 , \ldots,m_r)$ denote the linear system of surfaces of degree $d$
with multiplicities at least $m_i$ at $P_i$ for all $i=1,\ldots,r$ as well as the corresponding sheaf.

By abuse of notation, on $X_r$, $\cl_3(d ; m_1 , \ldots,m_r)$ also denotes the invertible sheaf
$\pi^*(\co_{\p^3}(d)) \otimes \co_{X_r}(-m_1E_1 - \cdots -m_rE_r)$ and the corresponding complete linear
system $|dE_0 - m_1E_1 - \cdots m_rE_r|$.

Analoguously, on $\p^3$, $\cl_3(d ; m_1^{n_1} , \ldots,m_t^{n_t})$ denotes the linear system
of surfaces of degree $d$ with multiplicities at least $m_i$ at $n_i$ of the points on $D_Q$
as well as the corresponding sheaf. Again, the same notation is used to denote
the associated complete linear system and invertible sheaf on $X_r$.

The virtual dimension of the linear system
$\cl = \cl_3(d ; m_1 , \ldots,m_r)$ on $\p^3$ as well as on $X_r$ is defined as
$$ \vdim(\cl) :=\binom{d+3}{3}-\sum_{i=1}^r\binom{m_i+2}{3} -1. $$
The expected dimension of $\cl$ is then given by
$$\edim(\cl) := \max \{ -1 , \vdim(\cl)\}.$$
It is then clear that $\dim(\cl) \geq \edim(\cl) \geq \vdim(\cl)$,
and the system $\cl$ is called {\em special} if $\dim(\cl) >
\edim(\cl)$. The system $\cl$ is associated to the sections of the
sheaf $\co_{\p^n}(d)\otimes\ci_Z$, where $Z=\sum m_ip_i$ is the
zero-dimensional scheme of fat points. From the cohomology exact
sequence associated to
\[
\xymatrix{ 0 \ar[r] & \cl\ar[r] & \co_{\p^n}(d)\ar[r] & \co_Z
\ar[r] & 0,}
\]
we obtain that $h^i(\cl) = 0$ for $i=2,3$. Therefore $v(\cl) =
h^0(\cl)-h^1(\cl)-1$, so that a non-empty system is special if and
only if $h^1(\cl)>0$.

Note that the strict transform $Q_r$ of ${\bar{Q}}$ on $X_r$
is a divisor of $\cl_3(2 ; 1^r)$, and $Q_r$ is just the blowing-up of $\bar{Q}$ along $P_1,\ldots,P_r$.
So $\Pic(Q_r) = \langle f_1, f_2, e_1, \ldots, e_r \rangle$,
with $f_1$ and $f_2$ the pullbacks of the two rulings on $\bar{Q}$
and $e_1, \ldots, e_r$ the exceptional curves.
By $\cl_{Q_r}(a, b ; m_1 , \ldots ,m_r)$ we denote
the complete linear system $|a f_1 + b f_2 - m_1 e_1 - \ldots - m_r e_r |$,
and, as before, if some of the multiplicities are the same, we also use
the notation $\cl_{Q_r}(a, b ; m_1^{n_1} , \ldots ,m_t^{n_t})$.

Let $B_s$ be the blowing-up of $\p^2$ along $s$ general points of a smooth irreducible cubic,
then $\Pic B_s = \langle h , e'_1 , \ldots, e'_s \rangle$,
with $h$ the pullback of a line and $e'_l$ the exceptional curves.
By $\cl_2( d ; m_1 , \ldots ,m_s)$ we denote
the complete linear system $|d h - m_1 e'_1 - \ldots - m_s e'_s |$.
And again, as before, if some of the multiplicities are the same, we also use
the notation $\cl_{2}(d ; m_1^{n_t} , \ldots ,m_r^{n_t})$.
Note that $-K_{B_s} = \cl_2( 3 ; 1^s)$ contains a smooth irreducible divisor
which we will denote by $D_{B_s}$.

On $B_s$, a system $\cl_2( d ; m_1 , \ldots ,m_s)$ is said to be in standard form if
$d \geq m_1 + m_2 + m_3$ and $m_1 \geq m_2 \geq \cdots \geq m_s \geq 0$;
and it is called standard it there exists a base
$\langle \hti , \et_1 , \ldots, \et_s \rangle$
of $\Pic B_s$
such that
$\cl_2( d ; m_1 , \ldots ,m_s) = |\tilde{d} \hti - \tilde{m}_1 \et_1 - \ldots - \tilde{m}_s \et_s |$
is in standard form.

As explained in~\cite[\S 6]{DL1}, the blowing-up $Q$ of the quadric along 1 general point
can also be seen as a blowing-up of the projective plane
along $2$ general points,
and
$$\cl_{Q}(a, b ; m) = \cl_2( a+b-m ; a-m , b-m ).$$
So, in particular $-K_Q = \cl_{Q}(2, 2 ; 1) = \cl_2( 3 ; 1^2 ) = -K_{B_2}$.
Obviously, this means that our blowing-up $Q_r$ can also be seen as a $B_{r+1}$
and
$$\cl_{Q_r}(a, b ; m_1, m_2, \ldots, m_r) = \cl_2( a+b-m_1 ; a-m_1 , b-m_1, m_2, \ldots, m_r ).$$
This implies in particular that we can apply the results from \cite{BH1} and \cite{BH2}.

\section{Non-speciality}

\begin{theor}\label{ns}
Consider $\cl = \cl_3(d ; m_1 , \ldots,m_r)$ on $X_r$ with $2d \geq m_1 + m_2 + m_3 + m_4$
and $m_1 \geq m_2 \geq \cdots \geq m_r \geq 0$.
Then $h^1(\cl) = 0$ \snvs if
\begin{itemize}
\item[(1)] $d \geq m_1 + m_2 -1$ and
\item[(2)] $4d \geq m_1 + \cdots m_r$ if $r \geq 9$.
\end{itemize}
\end{theor}

\begin{proof}
We will assume that $m_r > 0$, since otherwise we can work on $X_{r'}$ with $r' := \max \{ i : m_i >0 \}$.
If $r \leq 8$ then the points $P_1,\ldots,P_r$ are general points of $\p^3$
and the statement follows from \cite[Theorem 5.3]{DL1}.
So we assume that $r \geq 9$
and consider the following exact sequence
$$
\xymatrix@1{
0 \ar[r] & \cl_3(d-2;m_1-1,\ldots,m_r-1) \ar[r] & \cl \ar[r] &
    \cl \otimes \co_{Q_r} \ar[r]  & 0
}
$$
Then $\cl \otimes \co_{Q_r} = \cl_{Q_r}(d,d;m_1,\ldots,m_r) =
\cl_2(2d-m_1 ; (d-m_1)^2 , m_2 ,\ldots, m_r)$. Proceeding as in
the proof of \cite[Lemma 5.2]{DL1} it is easily seen that this is
a standard class. Moreover $\cl_2(2d-m_1 ; (d-m_1)^2 , m_2
,\ldots, m_r) . K_{B_{r+1}} = -4d +m_1 + \cdots +m_r \leq 0$, so
we can apply \cite[Theorem 1.1 and Proposition 1.2]{BH1} to obtain
$h^1(\cl \otimes \co_{Q_r}) = 0$. On the other hand, one can
easily check that $\cl' := \cl_3(d-2;m_1-1,\ldots,m_r-1)$ still
satisfies the conditions of the theorem. Continuing like this
until the residue class $\cl' = \cl_3(d';m'_1,\ldots,m'_{r'})$ is
such that $r' \leq 8$. For this class we then know that $h^1(\cl')
= 0$ which gives us that $h^1(\cl)=0$.
\end{proof}

\section{Base point freeness}

\begin{theor}\label{bpf}
Consider $\cl = \cl_3(d ; m_1 , \ldots,m_r)$ on $X_r$ with $m_1 \geq m_2 \geq \cdots \geq m_r$.
Then $\cl$ is base point free on $X_r$ if and only if
the following conditions are satisfied \snvs
\begin{itemize}
\item[(1)] $m_r \geq 0$,
\item[(2)] $d \geq m_1 + m_2$ and
\item[(2)] $4d \geq m_1 + \cdots m_r +2$ if $r \geq 8$.
\end{itemize}
\end{theor}

\begin{proof}
First of all let us prove that the conditions are necessary.
Obviously, if $m_r < 0$ then $m_r E_r \subset \bs(\cl)$;
and if $d < m_1 + m_2$ then the strict transform of the line through $P_1$ and $P_2$
is contained in $\bs(\cl)$.
Now assume that (1) and (2) are satisfied, but $4d \leq m_1 + \cdots m_r +1$
and consider
\begin{equation}\label{seq}
\xymatrix@1{
0 \ar[r] & \cl_3(d-2;m_1-1,\ldots,m_r-1) \ar[r] & \cl \ar[r] &
    \cl \otimes \co_{Q_r} \ar[r]  & 0
}
\end{equation}
Proceeding as before, one can check that $\cl \otimes \co_{Q_r} =
\cl_2(2d-m_1 ;\allowbreak (d-m_1)^2 ,\allowbreak m_2 ,\ldots,
m_r)$ is standard and $\cl_2(2d-m_1 ; (d-m_1)^2 , m_2 ,\ldots,
m_r) . (-K_{B_{r+1}}) = 4d -m_1 - \cdots -m_r \leq 1$. Using the
results from \cite{BH1}, we obtain that $\cl \otimes \co_{Q_r}$
has base points, which will also be base points of $\cl$.

Now assume that all three conditions are satisfied (as before, we may even assume $m_r > 0$).
If $r \leq 8$, the result follows from
\cite[Theorem 6.2]{DL2}, so we assume that $r \geq 9$
and that the result holds for $r' < r$
Consider the exact sequence (\ref{seq}).
As before, one can see that $\cl \otimes \co_{Q_r} = \cl_2(2d-m_1 ; (d-m_1)^2 , m_2 ,\ldots, m_r)$
is standard and since
$\cl_2(2d-m_1 ; (d-m_1)^2 , m_2 ,\ldots, m_r) . K_{B_{r+1}} \leq -2$
we know that $\cl \otimes \co_{Q_r}$ is base point free (see \cite[Lemma 3.3(2)]{BH1}).
Also, because of theorem~\ref{ns}, we know that $h^1(\cl_3(d-2;m_1-1,\ldots,m_r-1))=0$,
so $\cl$ induces the complete linear system $\cl \otimes \co_{Q_r}$ on $Q_r$,
which means in particular that $\cl$ has no base points on $Q_r$.
On the other hand it is easily checked that $\cl_3(d-2;m_1-1,\ldots,m_r-1)$  still satisfies
the conditions of the theorem.
Continue like this untill you have $r' < r$ for the residue class (i.e. repeat this reasoning $m_r$ times).
Denote $\cl_3(d-2m_r;m_1-m_r,\ldots,m_{r-1}-m_r)$ by $\cl'$.
We then know that $\cl' + mQ_r \subset \cl$. Since $\cl$ has no base points on $Q_r$, and since
$\cl'$ is base point free by our induction hypothesis, we obtain that $\bs(\cl) = \emptyset$.
\end{proof}

\section{Very ampleness}

\begin{theor}\label{va}
Consider $\cl = \cl_3(d ; m_1 , \ldots,m_r)$ on $X_r$ with $m_1 \geq m_2 \geq \cdots \geq m_r$.
Then $\cl$ is very ample on $X_r$ if and only if
the following conditions are satisfied \snvs
\begin{itemize}
\item[(1)] $m_r > 0$,
\item[(2)] $d \geq m_1 + m_2 + 1$ ($d \geq m_1 +1$ if $r=1$; $d \geq 1$ if $r=0$) and
\item[(2)] $4d \geq m_1 + \cdots m_r + 3$ if $r \geq 9$.
\end{itemize}
\end{theor}

\begin{proof}
First of all let us note that the conditions are necessary.
Obviously, if $m_r \leq 0$ then $\cl$ cannot separate on $E_r$;
and if $d \leq m_1 + m_2$ then $\cl$ cannot separate
a zero-dimensional subscheme $Z$ of length 2 of
the strict transform of the line through $P_1$ and $P_2$.
In case $4d \leq m_1 + \cdots m_r + 2$, one can see that $\cl$ cannot separate
$Z$ if it is contained in $D_{Q_r}$.

Now assume that all three conditions are satisfied.

First of all consider the exact sequence
\[
\xymatrix@1{
0 \ar[r] & \cl_3(d;m_1+1, m_2,\ldots ,m_r) \ar[r] & \cl \ar[r] &
    \cl \otimes \co_{E_1} \ar[r]  & 0
}
\]
Because of Theorem~\ref{ns} we know that $h^1(\cl_3(d;m_1+1, m_2,\ldots ,m_r))=0$,
so $\cl$ induces the complete linear system $\cl \otimes \co_{E_1}$ on $E_1$.
Since $\cl \otimes \co_{E_1} = \cl_2(m_1)$, we see that $\cl$ separates on $E_1$.
Naturally, a similar reasoning can be done for any $E_i$, so $\cl$ separates on every $E_i$.

Moreover, one can easily check that $\cl(E_i) := \cl_3(d ;
m_1,\allowbreak \ldots, m_{i-1} ,\allowbreak m_i+1,\allowbreak
m_{i+1} ,\allowbreak \ldots,\allowbreak m_r)$ satisfies all the
conditions of Theorem~\ref{bpf}, so that $\cl(E_i)$ is base point
free on $X_r$. This means that $\cl$ can separate $Z$ if $\exists
i : Z \cap E_i \neq 0$ but $Z \not\subset E_i$.

Combining the previous two results, we see that we now only need to show that $\cl$ separates $Z$
with $Z \cap E_i = \emptyset$ for all $i=1,\ldots,r$.

In case $r =0,1$ or $2$, this is trivial.

Now let us assume that $r \geq 3$
and that the statement holds for $r' < r$.

First look at the case where $m_r = 1$
and consider the exact sequence
\[
\xymatrix@1{
0 \ar[r] & \cl_3(d-2;m_1-1,\ldots,m_{r-1}-1,0) \ar[r] & \cl \ar[r] &
    \cl \otimes \co_{Q_r} \ar[r]  & 0
}
\]
Proceeding as before, one can easily see that
$\cl \otimes \co_{Q_r}$ is standard.
Moreover, $4d - m_1 - \cdots - m_r \geq 3$ (if $r \geq 9$ this is condition (3) and if $3 \leq r \leq 8$
this follows from (2)), so \cite[Theorem 2.1]{BH2} implies that
$\cl \otimes \co_{Q_r}$ is very ample.
Since
$h^1(\cl_3(d-2;m_1-1,\ldots,m_{r-1}-1))=0$ (because of Theorem~\ref{ns}) we then obtain that
$\cl$ separates on $Q_r$.
Using Theorem~\ref{bpf} we also obtain that $\cl' := \cl_3(d-2;m_1-1,\ldots,m_{r-1}-1, 0)$
is base point free, which implies that $\cl$ separates $Z$ if $Z \cap Q_r \neq \emptyset$.
Let $r'  = \max \{ i : m_i >1\}$ (or $r' = 0$ if all $m_i =1$),
then one can easily check that
$\cl'$ satisfies all the conditions of the theorem on $X_{r'}$.
So, using the induction hypothesis, we have that $\cl'$ is very ample on $X_{r'}$.
But since $Z$ on $X_r$ is disjoint with all $E_i$, $Z$ corresponds with a zero-dimensional subscheme
on $X_{r'}$ (also disjoint with all $E_i$).
So may may conclude that $\cl$ separates any $Z$.

Now we assume $m_r > 1$ and we assume that the
statement holds for $m'_r < m_r$.
Consider the exact sequence
\[
\xymatrix@1{
0 \ar[r] & \cl_3(d-2;m_1-1,\ldots,m_{r}-1) \ar[r] & \cl \ar[r] &
    \cl \otimes \co_{Q_r} \ar[r]  & 0
}
\]
Proceeding similarly as for the case $m_r=1$ one can easily see that $\cl$ separates all $Z$.
\end{proof}

\begin{rem}
One can also use Theorem~\ref{bpf} and \cite[Theorem~2.1]{BS} to obtain a very ampleness result.
However in this way only part of the complete class of very ample systems on $X_r$ are obtained.
\end{rem}

\begin{rem}
Let $A_1,\ldots,A_r$ be general points on $\p^3$, let $Y_r$ be the blowing-up of
$\p^3$ along those $r$ general points and let $\cl_3(d ; m_1,\ldots,m_r)$
($m_1 \geq m_2 \geq \cdots \geq m_r$)
denote the complete linear system $|dE_0 - m_1 E_1 - \cdots -m_r E_r|$ on $Y_r$.
Since the very ampleness is an open property, Theorem~\ref{va} implies that
$\cl_3(d ; m_1,\ldots,m_r)$ is very ample on $Y_r$ if
$m_r > 0$,
$d \geq m_1 + m_2 + 1$ ($d \geq m_1 +1$ if $r=1$; $d \geq 1$ if $r=0$) and
$4d \geq m_1 + \cdots m_r + 3$ if $r \geq 9$.
Of course the third condition will now no longer be a necessary condition.
\end{rem}

%\bibliographystyle{alpha}

%\bibliography{BaseLoc-P3.bib}

\end{document}